\documentclass[12pt]{amsart}

\usepackage{amsthm,amsmath,amsfonts}
\usepackage[all]{xy}
\usepackage{psfrag}
\usepackage[dvips]{graphicx}

\usepackage[latin1]{inputenc}

\makeatletter
\def\th@alexnormal{%
\let\thm@indent\noindent 
\thm@headfont{\bfseries}

\normalfont
}
\makeatother

\makeatletter
\def\th@alexit{%
\let\thm@indent\noindent 
\thm@headfont{\bfseries}
\normalfont
\fontshape{it}
\selectfont
}
\makeatother

\theoremstyle{alexit}
\newtheorem{theorem}{Theorem}[subsection]
\newtheorem{proposition}[theorem]{Proposition}
\newtheorem{corollary}[theorem]{Corollary}
\newtheorem{lemma}[theorem]{Lemma}
\newtheorem{conjecture}[theorem]{Conjecture}

\theoremstyle{alexnormal}
\newtheorem{definition}[theorem]{Definition}
\newtheorem{remark}[theorem]{Remark}

\numberwithin{equation}{subsection}

\begin{document}
\author{Alexandre Girouard}
\title[Concentration to points, conformal degeneration and
$\lambda_1$]
{Fundamental tone, concentration of density to points and conformal degeneration on surfaces}
\date{\today}

\address{D\'epartement de Math\'ematiques et
Statistique, Universit\'e de Montr\'eal, C. P. 6128,
Succ. Centre-ville, Montr\'eal, Canada H3C 3J7}
\email{girouard@dms.umontreal.ca}

\maketitle

\begin{abstract}
  We study the effect of two types of degeneration of the Riemannian
  metric on the first eigenvalue of the Laplace operator on
  surfaces. In both cases we prove that the first eigenvalue of the
  round sphere is an optimal asymptotic upper bound. The first type of
  degeneration is concentration of the density to a point within a
  conformal class. The second is degeneration of the 
  conformal class to the boundary of the moduli space on the torus and
  on the Klein bottle. In the latter, we follow the outline proposed
  by N. Nadirashvili in 1996.
\end{abstract}

\section{Introduction}
Given a Riemannian metric $g$ on a closed surface $\Sigma$, 
let the spectrum of the Laplace operator
$\Delta_g$ acting on smooth functions be the sequence
$$0=\lambda_0(g)<\lambda_1(g)\leq\lambda_2(g)\leq\cdots\leq\lambda_k(g)\leq\cdots\nearrow\infty$$
where each eigenvalue is repeated according to its multiplicity. The
first nonzero eigenvalue  $\lambda_1(g)$ is called the
\emph{fundamental tone} of $(\Sigma,g)$.
Let $\mathcal{R}(\Sigma)$ be the space of Riemannian metrics on
$\Sigma$ with total area one. 
We are interested in the asymptotic behavior of the functional
$$\lambda_1:\mathcal{R}(\Sigma)\rightarrow ]0,\infty[$$
under two types of degeneration of sequences of metrics of area one
$(g_n)\subset\mathcal{R}(\Sigma)$ derscribed below.

\subsection{Concentration to points}
It is expected that a metric maximizing
$\lambda_1:\mathcal{R}(\Sigma)\rightarrow]0,\infty[$ 
has a lot of symmetries. On the sphere, the torus and
the projective plane for example, the $\lambda_1$-maximizing metrics
are the standard homogeneous ones. We consider the opposite situation
where the distribution of mass of a sequence of metrics concentrate to
a point, developping a $\delta$-like singularity.
\begin{definition}
  A sequence $(g_n)\subset \mathcal{R}(\Sigma)$  is said to
  \emph{concentrate to the point $p\in \Sigma$} if for each
  neighborhood $\mathcal{O}$ of $p$
  $$\lim_{n\rightarrow\infty}\int_{\mathcal{O}}dg_n=1.$$  
\end{definition}

\noindent
\textbf{Question}. Does concentration to a point impose any restriction
to the asymptotic behavior of the eigenvalues of the Laplace operator
$\Delta_{g_n}$ on the surface $\Sigma$?\\

Without any further constraints, the answer is no.
\begin{proposition}\label{no}
  For any metric $g_0$, and any point $p\in\Sigma$ there exists
  a sequence $(g_n)$ of pairwise isometric metrics
  concentrating to $p$. In particular the metrics
  $(g_n)$ are isospectral.
\end{proposition}
Under the additional assumption that the metrics $g_n$ are conformally
equivalent, we obtain an optimal asymptotic upper bound on the
fundamental tone.
\begin{theorem}\label{punktthm}
  Let $[g]=\{\alpha g\ |\ \alpha\in C^\infty(\Sigma), \alpha>0\}$ be a
  conformal class on a closed surface $\Sigma$.
  \begin{itemize}
  \item[a)] For any sequence $(g_n)$ of metrics in the conformal class $[g]$ which concentrates to a point $p\in\Sigma$,
    $$\limsup_{n\rightarrow\infty}\lambda_1(g_n)\leq 8\pi.$$
  \item[b)] For any point $p\in\Sigma$, there exists a sequence 
    $(g_n)$ of metrics of area one in the conformal class $[g]$
    concentrating to $p$ and such 
    that $$\lim_{n\rightarrow\infty}\lambda_1(g_n)=8\pi.$$
  \end{itemize}
\end{theorem}
Proposition~\ref{no} and Theorem~\ref{punktthm} will be proved in
section~\ref{sectionPunkt}.

\subsection{Conformal degeneration}
Given a conformal class $[g]$ on the torus $T^2$, define 
$$\nu([g]):=\sup_{\tilde{g}\in\mathcal{R}(T^2)\cap [g]}\lambda_1(\tilde{g}).$$
This corresponds to the first conformal eigenvalue of Colbois and El
Soufi~\cite{colbois:2}.
Let
$$\mathcal{M}:=\{a+ib\in\mathbb{C}\ |\ 0\leq a\leq 1/2, a^2+b^2\geq 1, b>0\}.$$
Any metric on $T^2$ is
conformally equivalent to a flat torus $\mathbb{C}/\Gamma$ for some
lattice $\Gamma$ of $\mathbb{C}$ generated by $1\in\mathbb{C}$ and
$a+ib\in\mathcal{M}$.
\begin{figure}[h]
  \centering
  \psfrag{1}[][][1]{$1$}
  \psfrag{a+ib}[][][1]{$a+ib$}
  \psfrag{M}[][][1]{$\mathcal{M}$}
  \psfrag{*}[][][1]{$\bf{*}$}
  \includegraphics[width=6cm]{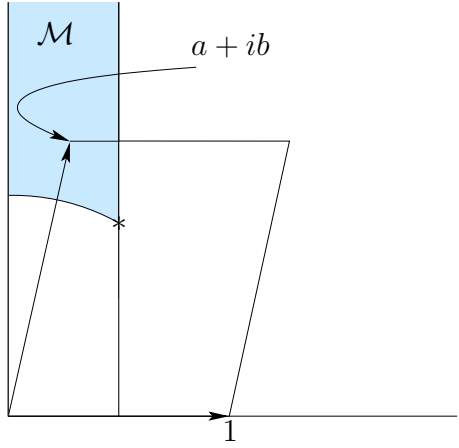}
  \caption{Moduli space of conformal structures on $T^2$}
  \label{figmodtorus}
\end{figure}
It follows that $\mathcal{M}$ is a natural representation of the
moduli space $\mathcal{M}(T^2)$ of conformal classes on the torus (See
Figure~\ref{figmodtorus}). 
\begin{definition}
  A sequence of metrics on the torus $T^2$ is \emph{degenerate} if the
  corresponding sequence $(a_n+ib_n)\subset\mathcal{M}$ satisfies
  $\lim_{n\rightarrow\infty}b_n=\infty.$
\end{definition}
\begin{theorem}\label{monthmNad}
  If a sequence $(g_n)$ of Riemannian metrics of area one on the torus
  is degenerate, then 
  $$\lim_{n\rightarrow\infty}\nu([g_n])=8\pi.$$
  In particular, $$\limsup_{n\rightarrow\infty}\lambda_1(g_n)\leq 8\pi.$$
\end{theorem}
The proof will be presented in section \ref{preuveconftorus}.
We complete the outline proposed by Nadirashvili in \cite{nad:1}.
In particular, a detailed version of a concentration lemma
(Lemma \ref{lemmeConcentration}) is used and a Dirichlet energy
estimate on long cylinders is proposed (Lemma~\ref{nrjboundary}
and Lemma~\ref{nrjdirich}).\\

\subsection{Maximization of $\lambda_1$ on surfaces}
One motivation for Theorem~\ref{monthmNad} is the role that it plays
in $\lambda_1$-maximization on the torus.
More generally, it is natural to ask which metric on a closed surface
$\Sigma$, if any, maximizes the fundamental tone $\lambda_1$ in the
space $\mathcal{R}(\Sigma)$ of Riemannian metrics of area one.
It is known that
$$\Lambda(\Sigma):=\sup_{g\in\mathcal{R}(\Sigma)}\lambda_1(g)$$
is finite \cite{hersch:1}, \cite{yau:1}, \cite{li:1} for any closed
surface $\Sigma$, but the explicit value of $\Lambda(\Sigma)$ has only
been computed for few surfaces.

The first such result was obtained by
Hersch~\cite{hersch:1} in 1970. He proved that the round metric
$g_{S^2}$ of area one is the unique maximum of $\lambda_1$ on
$\mathcal{R}(S^2)$. The proof rests on the fact that by Riemann's
uniformization theorem  any two metrics on the sphere are conformally
equivalent.

In 1973, Berger~\cite{berger:2} proved that among flat metrics on the
torus, $\lambda_1$ is maximized by the flat equilateral metric
$g_{eq}$, that is, the metric induced from the quotient of
$\mathbb{C}$ by the lattice generated by $1$ and
$e^{i\pi/3}$(indicated by the $\bf{*}$ in
Figure~\ref{figmodtorus}). He conjectured that this metric is a global
maximum of $\lambda_1$ over all Riemannian metrics of total area
one. In 1996, Nadirashvili \cite{nad:1} proposed a method of proof.\\

\noindent
\textbf{Nadirashvili's approach.}
Start with a maximising sequence $(g_n)$,
i.e. such that $\lambda_1(g_n)\rightarrow\Lambda(T^2)$,
and show that it admits a subsequence converging to a real analytic
metric $\overline{g}$. Nadirashvili~\cite{nad:1} has proved that a
metric maximizing $\lambda_1$ on a surface $\Sigma$ is also
$\lambda_1$-minimal, which means that 
$(\Sigma,\overline{g})$ is minimally immersed
in a round sphere by its first eigenfunctions.
Montiel and Ros~\cite{montiel} have proved that for any surface
different from the sphere (in particular for the torus!), the
isometry group of a $\lambda_1$-minimal metric $\overline{g}$ coincide
with its group of conformal transformations. See also
\cite{soufi:1}. Since the group of conformal transformations of a
torus acts transitively, this implies that the curvature of
$\overline{g}$ is constant and therefore zero by Gauss-Bonnet
theorem. The above result of Berger \cite{berger:2} completes the
proof.

The first step in showing that $(g_n)$ admits a convergent subsequence
is to prove that the associated sequence $([g_n])$ of conformal
classes admits a converging subsequence. Since 
$\lambda_1(g_{eq})=8\pi^2/\sqrt{3}>8\pi$, Theorem~\ref{monthmNad}
implies that the corresponding sequence $(a_n+ib_n)$ is bounded and
must therefore admits a convergent subsequence.

\begin{remark}
  \begin{itemize}
  \item[]
  \item[-]   Explicit $\lambda_1$-maximal metrics are
    also known for the projective space~\cite{li:1} and the Klein
    bottle~\cite{jakobson:2},~\cite{soufi:2}. There is a conjecture for
    surfaces of genus two~\cite{jakobson:1}.
  \item[-] The existence of analytic $\lambda_1$-maximal metrics has recently
    been  used in~\cite{jakobson:2} and~\cite{soufi:2}.
  \end{itemize}
\end{remark}

\subsection{Conformal degeneration on the Klein bottle}
Define two affine transformations $t_b$ and $\tau$ of $\mathbb{C}$ by
\begin{align*}
  t_b(x+iy)&=x+i(y+b),\\
  \tau(x+iy)&=x+\pi-iy.
\end{align*}
Let $G_b$ be the group of transformations generated by $t_b$ and $\tau$.
\begin{lemma}
  Any Riemannian metric $g$ on the Klein bottle $\mathbb{K}$ is
  conformally equivalent to one of the standard flat models
  $$K_b:=\mathbb{C}/G_b.$$
  In other words, there exist a smooth function
  $\alpha:K_b\rightarrow]0,\infty[$ such that $(\mathbb{K},g)$ 
  is isometric to $\left(K_b,\alpha(dx^2+dy^2)\right).$
\end{lemma}
It follows that the moduli space of conformal classes on the Klein
bottle is identified with the set of positive real numbers.
\begin{figure}[h]
  \centering
  \psfrag{1}[][][1]{$1$}
  \psfrag{b}[][][1]{$b$}
  \psfrag{p}[][][1]{$\pi$}
  \includegraphics[width=10cm]{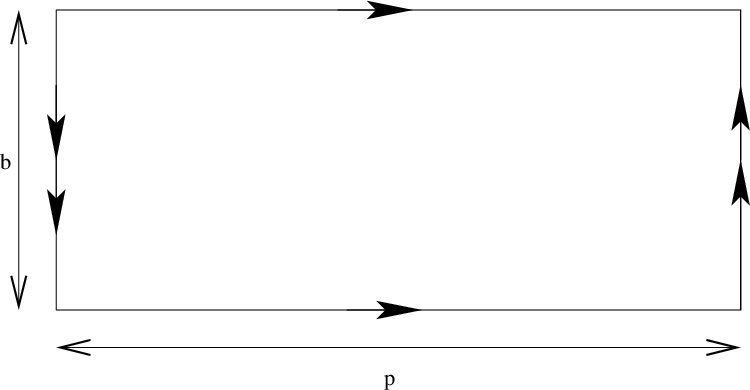}
\end{figure}
\begin{theorem}\label{kleindegeneration}
  Let $(g_n)\subset\mathcal{R}(\mathbb{K})$ be a
  sequence of metrics on the Klein bottle.
  \begin{itemize}
  \item[1.] If $\lim_{n\rightarrow\infty}b_n=0$, then
    $\limsup_{n\rightarrow\infty}\lambda_1(g_n)=8\pi.$
  \item[2.] If $\lim_{n\rightarrow\infty}b_n=\infty$, then
    $\limsup_{n\rightarrow\infty}\lambda_1(g_n)\leq 12\pi.$
  \end{itemize} 
\end{theorem}
The proof will be presented in section~\ref{sectionklein}, we follow
the outline proposed by Nadirashvili in \cite{nad:1} and in a private
communication. The first case is very similar to the corresponding
result for the torus (Theorem~\ref{monthmNad}). The second case uses
the fact that the standard metric on $\mathbb{R}P^2$ is minimally
embedded in $S^4$ by its first eigenfunctions to obtain an estimate
on the Dirichlet energy of a test function via a theorem of Li and Yau
\cite{li:1} on conformal area of minimal surfaces. 

\subsection{Friedlander and Nadirashvili  invariant}
For a closed manifold $M$ of dimension at least $3$, Colbois and Dodziuk
\cite{colbois:1} proved that the first eigenvalue is unbounded on the
set of Riemannian metrics of area one, that is
$\Lambda(M)=+\infty$. On the other hand, it is known that the supremum
$\nu(C)$ of $\lambda_1$ restricted to metrics of area one in any fixed
conformal class $C$ is finite \cite{soufi:3}. Friedlander and
Nadirashvili~\cite{nad:2} introduced the differential invariant
$$I(M):=\inf\left\{\nu(C)\ \Bigl|\Bigr.\ C \mbox{ is a conformal class
    on } M\right\}$$
and proved that it satisfies
$$I(M)\geq \lambda_1(S^n,g_{S^n})$$
where $g_{S^n}$ is the round metric of area one on the sphere $S^n$.
It would be interesting to know if this invariant distinguishes
nonequivalent differential structures. In fact, it is very difficult
to compute $I(\Sigma)$ explicitely, even for surfaces. Since any two
metrics on $S^2$ are conformally equivalent, it is obvious that
$I(S^2)=8\pi$. For similar reasons, $I(\mathbb{R}P^2)=12\pi$. The
invariant for the torus and for the Klein bottle are obtained as
corollaries to Theorem~\ref{monthmNad} and
Theorem~~\ref{kleindegeneration}.
\begin{corollary}
  The Friedlander-Nadirashvili invariants  the torus and the Klein
  bottle are $8\pi$.
\end{corollary}
Let us also mention the following conjecture.
\begin{conjecture}[Friedlander-Nadirashvili]
  For any closed  surface $\Sigma$ other than the projective plane
  $\mathbb{R}P^2$, $I(\Sigma)=8\pi.$
\end{conjecture}

\section{Analytic background}
Let $\Sigma$ be a closed surface.
In dimension two the Dirichlet energy of a function $u\in C^\infty(\Sigma)$
$$D(u)=\int_\Sigma |\nabla_gu|^2\,dg$$
is invariant under conformal diffeomorphisms. In order to estimate the
first eigenvalue of the Laplace operator $\Delta_g$,  the following
variational characterization will be used:
\begin{gather}\label{varcharaclambda1}
  \lambda_1(g)=\inf
  \left\{\frac{D(u)}{\int_\Sigma u^2\,dg}\
    \Bigl|\Bigr.\ u\in C^\infty(\Sigma), u\neq 0, 
    \int_\Sigma u\,dg=0\right\}.
\end{gather}

\subsection{Dirichlet energy estimate on thin cylinders}
The main technical tool that we use is an  estimate on the
Dirichlet energy of harmonic functions on long cylinders in terms of
their restrictions to the boundary circles.
\begin{lemma}\label{nrjdirich}
  Let $\Omega=(0,L)\times S^1$ with
  $S^1=\mathbb{R}/2\pi\mathbb{Z}$. Consider
  $f\in C^{\infty}(\overline{\Omega})$ such that $f(0,\theta)=0$ and
  $|f(L,\theta)|\leq 1$.
  Let $u(\theta)=f(L,\theta)$. If $f$ is harmonic, then
  $$\int_\Omega|\nabla f|^2\leq\frac{2\pi}{L}+\coth(L)\int_{\theta=0}^{2\pi}|u'(\theta)|^2d\theta.$$
\end{lemma}
The idea is to express the Fourier series of the function $f$ in terms of
the Fourier series of $u$. This is similar to Hurwitz's proof of the
isoperimetric inequality \cite{hurwitz:1}.
\begin{proof}
  Let the Fourier representation of $u$ be
  $$u(\theta)=K+\sum^\infty_{k=1}\bigl(a_k\cos(k\theta)+b_k\sin(k\theta)\bigr).$$
  Direct computation shows that $f$ admits the following representation:
  $$f(x,\theta)=\frac{Kx}{L}+\sum^\infty_{k=1}\frac{\sinh(kx)}{\sinh(kL)}\Bigl(a_k\cos(k\theta)+b_k\sin(k\theta)\Bigr)$$
  and integration by parts leads to
  \begin{align*}
    \int_\Omega|\nabla f|^2&=-\int_\Omega f\Delta
    f+\int_{\partial\Omega}f\frac{\partial f}{\partial\nu}\\
    &=\int_{\theta=0}^{2\pi}u(\theta)\partial_xf(x,\theta)d\theta\Bigl|\Bigr._{x=L}\\
    &=\int_{\theta=0}^{2\pi}\Bigl(K+\sum^\infty_{k=1}a_k\cos(k\theta)+b_k\sin(k\theta)\Bigr)\times\\
    &\ \ \ \ \ \ \ \ \ \ 
    \Bigl(\frac{K}{L}+\sum^\infty_{k=1}k\coth(kL)\left(a_k\cos(k\theta)+b_k\sin(k\theta)\right)\Bigr)d\theta\\
    &=2\pi\frac{K^2}{L}+\pi\sum_{k=1}^\infty k\coth(kL)(a_k^2+b_k^2)
  \end{align*}
  For $x>0$,
  $$\frac{d}{dx}\coth(x)=-\frac{4}{(e^x-e^{-x})^2}<0$$
  so that for each $k\geq 1$,
  $$\coth(kL)\leq\coth(L).$$
  It follows that
  \begin{align*}
    \int_\Omega|\nabla f|^2&\leq 2\pi\frac{K^2}{L}+
    \pi\coth(L)\sum_{k=1}^\infty k(a_k^2+b_k^2)\\
    &\leq 2\pi\frac{K^2}{L}+\pi\coth(L)\sum_{k=1}^\infty k^2(a_k^2+b_k^2)\\
    &=2\pi\frac{K^2}{L}+\coth(L)\int_{\theta=0}^{2\pi}|u'(\theta)|^2d\theta
  \end{align*}
  where $|K|=\frac{1}{2\pi}|\int_{\theta=0}^{2\pi}u(\theta)d\theta|\leq 1$ since
  $u(\theta)\in[-1,1]$.
\end{proof}
A simple conformal change of coordinates is used to extend
Lemma~\ref{nrjdirich} to the situation where the boundary circle has
arbitrary length.
\begin{corollary}\label{coro_nrjdirich}
  Suppose the hypothesis of Lemma~\ref{nrjdirich} holds with the
  circle $\mathbb{R}/2\pi\mathbb{Z}$ replaced by
  $\mathbb{R}/r\mathbb{Z}$. For any $L>0$  
  \begin{gather}\label{ineqdirichlet}
    \int_\Omega|\nabla f|^2\leq \frac{r}{L}+
    \frac{r}{2\pi}\coth(\frac{2\pi L}{r})\int_{x=0}^{r}|u'(x)|^2dx.
  \end{gather}
\end{corollary}

\subsection{Conformal renormalization of centers of mass}
It is possible to conformally move any nonsingular distribution of
mass on the sphere $S^n\subset\mathbb{R}^{n+1}$ in such a way that
its center of mass becomes the origin of $\mathbb{R}^{n+1}$.
\begin{lemma}[Hersch Lemma]\label{lemmealaHerschgeneral}
  Let $\mu$ be a measure on the sphere $S^n$. If the support of $\mu$
  is not a point, then there exists a conformal transformation $\tau$ of
  the sphere $S^n$ such that
  \begin{gather}\label{gravityOrigin}
    \int_{S^n}\pi\circ\tau\ d\mu=0
  \end{gather}
  where $S^n\overset{\pi}{\hookrightarrow}\mathbb{R}^{n+1}$ is the
  standard embedding.
\end{lemma}
This lemma was obtained by Hersch~\cite{hersch:1} in 1970, see also
\cite{schoen:1}. It is proved using a topological argument similar to
the proof of Brouwer fixed point theorem.
\begin{corollary}\label{coroAlaHerschplongement}
  Let $\mu$ be a measure on a surface $\Sigma$. Consider an
  embedding $\phi:\Sigma\rightarrow S^n$. If the support of $\mu$
  is not a point, then there
  exists a conformal transformation $\tau$ of $S^n$ such that
  \begin{gather*}
    \int_{\Sigma}\pi\circ\tau\circ\phi\,d\mu=0.
\end{gather*}  
\end{corollary}
\begin{proof}
  The result follows from application of Hersch Lemma to the
  push-forward measure $\phi_*\mu$ since
  $\int_{S^n}f\,d(\phi_*\mu)=\int_{\Sigma}f\circ\phi\,d\mu$
  for any smooth function $f$.
\end{proof}

\section{Conformal degeneration on the torus}\label{preuveconftorus}
The goal of this section is to prove Theorem~\ref{monthmNad}. 

\subsection{Moduli space of tori} 
For any $b>0$, let 
$$T_b:=\left\{[x+iy]\in\mathbb{C}/\mathbb{Z}\ \Bigl|\Bigr.\
  -b/2<y<b/2\right\}$$
be a cylinder of length $b$. Given $a+ib\in\mathcal{M}$, let
$\Gamma_{a,b}$ be the lattice of $\mathbb{C}$ generated by~1 and
$a+ib$. Consider the group $G_{a,b}$ of transformations of $T_\infty$
generated by
$$[x+iy]\mapsto[x+a+i(y+b)].$$
The cylinder $T_b$ is a fundamental domain of this
action and the torus $\mathbb{C}/\Gamma_{a,b}$ can also be obtained as
$T_\infty/G_{a,b}$.

\subsection{Concentration on thin cylinders}
In order to make notation less cumbersome, consider a sequence
$(g_n)$ of metrics
such that $b_n=n$. The first eigenvalue of the flat torus
corresponding to $g_n$,
$$\lambda_1(T_n/G_{a_n,n})=\frac{4\pi^2}{n^2}$$
tends to zero with $n$ going to
infinity~\cite{berger:1}. Imposing a uniform positive lower bound $\lambda_1(g_n)\geq
K>0$ on the first eigenvalues for $g_n$ should therefore imply that
$g_n$ is ``far from being flat''. The next lemma makes this intuitive
idea precise by showing that the Riemannian measures $dg_n$ concentrate on
relatively thin cylindrical parts in $(T^2,g_n)$.
\begin{lemma}\label{lemmeConcentration}
  If there exists $K>0$ such that
  $$\liminf_{n\rightarrow\infty}\lambda_1(g_n)\geq K,$$
  then for any $\epsilon>0$, there exists $N\in\mathbb{N}$ such that for
  any $n\geq N$, 
  $$\max\bigl\{\int_{T_{3n/4}}dg_n,\int_{T_n\setminus T_{n/4}}dg_n\bigr\}\geq 1-\epsilon.$$
\end{lemma}
Let $A_n^\epsilon$ be the maximizing cylinder: either $T_{3n/4}$ or
$T_n\setminus T_{n/4}$. This lemma says that most of the mass
(i.e. $1-\epsilon$) is concentrated on a cylinder whose lenth is $3/4$
of the total length. Without loss of generality, we will suppose
$A_n^\epsilon=T_{3n/4}$.
\begin{proof}
  The function
  $$\gamma_n([x+iy])=\cos(\frac{2\pi y}{n})$$
  is a first eigenfunction on the flat torus $T_n/G_{a_n,n}$
  corresponding to $g_n$.
    \begin{figure}[h]
    \centering
    \psfrag{c-d}[][][.7]{$c_n-\delta$}
    \psfrag{c+d}[][][.7]{$c_n+\delta$}
    \psfrag{A}[][][.7]{$P_n^\delta$}
    \psfrag{c}[][][.7]{$c$}
    \psfrag{c_n}[][][.7]{$c_n$}
    \includegraphics[width=10cm]{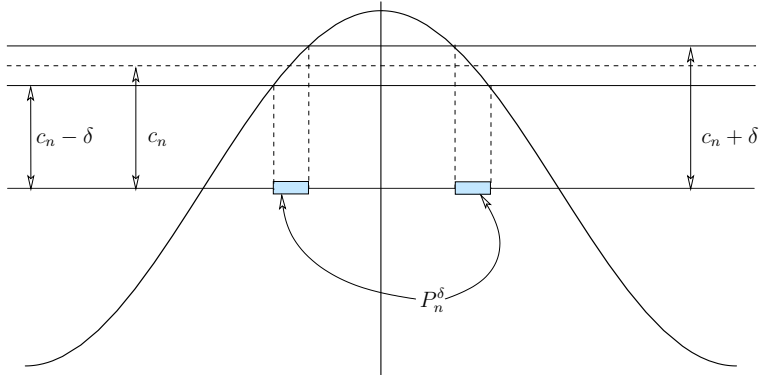}
    \caption{Concentration}
    \label{figconcentration}
  \end{figure}
  Let $$c_n=\int_{T_{n}}\gamma_n(x+iy)\,dg_n$$ and define
  $h_n:T_{n}\rightarrow\mathbb{R}$ by
  $$h_n=\gamma_n-c_n,$$
  which satisfies
  $$\int_{T_n}h_n\,dg_n=0.$$
  Given $\delta>0$, let
  $P_n^\delta=h_n^{-1}([-\delta,\delta])$
  and $Q_n^\delta=T_{n}\setminus P_n^\delta$. 
  Since $\delta^2\leq h_n^2$ on $Q_n^\delta$,
  $$\delta^2\int_{Q_n^\delta}dg_n\leq\int_{Q_n^\delta}h_n^2\,dg_n\leq\int_{T_{n}}h_n^2\,dg_n$$
  Using the 
  variational characterization of $\lambda_1(g_n)$ and the conformal
  invariance of the Dirichlet energy,
  \begin{align*}
    \delta^2\int_{Q_n^\delta}dg_n&\leq\frac{\int_{T_{n}}|\nabla h_n|^2\,dg_n}{\lambda_1(g_n)}
    \leq\frac{\int_{T_{n}}|\nabla \gamma_n|^2\,dg_n}{K}=\frac{4\pi^2}{Kn^2}\rightarrow 0.
  \end{align*}
  So that, for each $\delta>0$,
  $$\lim_{n\rightarrow\infty}\int_{Q_n^\delta}dg_n=0$$ and
  $$\lim_{n\rightarrow\infty}\int_{P_n^\delta}dg_n=1.$$
  Observe that
  $$P_n^\delta=\left\{[x+iy]\ \bigl|\bigr.\ c_n-\delta\leq cos(2\pi y/n)\leq c_n+\delta\right\}.$$
  For $\delta>0$ small enough, $P_n^{\delta}$ is the union of two
  intervals around $\pm\arccos(c_n)$ whose lengths
  are at most $\frac{n}{10}$.
  The set $P_n^\delta$ will either be included in $T_{3n/4}$ or in $
  T_{n}\setminus T_{n/4}$ (see Figure~\ref{figconcentration}).
\end{proof}

\subsection{Transplantation to the sphere}
Let $\sigma:\mathbb{C}\rightarrow S^2$
be the stereographic parametrization of the sphere by its equatorial plane
\begin{gather}\label{stereoparam}
  \sigma(u+iv)=\frac{1}{1+u^2+v^2}(2u,2v,u^2+v^2-1)
\end{gather}
and define the conformal equivalence
$\phi:\mathbb{C}/\mathbb{Z}\rightarrow S^2\setminus\{\mbox{poles}\}$ by
$$\phi([z])=\sigma(e^{-2\pi iz}).$$

\subsection{Renormalization of the centers of mass}
It follows from Corollary~\ref{coroAlaHerschplongement} that there exists
conformal transformations $\tau_n$ such that
\begin{gather}\label{renormcentermass}
  \int_{A_n^\epsilon}\pi\circ\tau_n\circ\phi\,dg_n=0
\end{gather}
where $S^2\overset{\pi}{\hookrightarrow}\mathbb{R}^3$ is the standard embedding.
For each $n\in\mathbb{N}$,
since $\pi_1^2+\pi_2^2+\pi_3^2=1$ on $S^2$, there exists an indice $i=i(n)\in\{1,2,3\}$
such that the function
\begin{gather}\label{deftorusun}
  u_n=\pi_i\circ\tau_n\circ\phi
\end{gather}
satisfies
\begin{gather}\label{tagada}
  \int_{A_n^\epsilon}u_n^2\,dg_n\geq\frac{1}{3}\int_{A_n^\epsilon}dg_n\geq\frac{1-\epsilon}{3}.
\end{gather}

\subsection{Test functions}
The function $u_n$ will be extended and perturbed to a function $f_n$
defined on $T_n$ and admissible for the variational
characterization~\eqref{varcharaclambda1} of $\lambda_1(g_n)$.

Let $I_n^-=[-7n/16,-6n/16]$ and $I_n^+=[6n/16,7n/16].$
\begin{figure}[h]
  \centering
  \psfrag{I_n^-}[][][.7]{$I_n^-$}
  \psfrag{I_n^+}[][][.7]{$I_n^+$}
  \psfrag{a_n^-}[][][.7]{$\alpha_n^-$}
  \psfrag{a_n^+}[][][.7]{$\alpha_n^+$}
  \psfrag{A_n}[][][.7]{$A_n^\epsilon$}
  \psfrag{3n/4}[][][.7]{$\frac{3}{4}n$}
  \psfrag{n/8}[][][.7]{$\frac{n}{8}$}
  \psfrag{n}[][][.7]{$n$}
  \includegraphics[width=10cm]{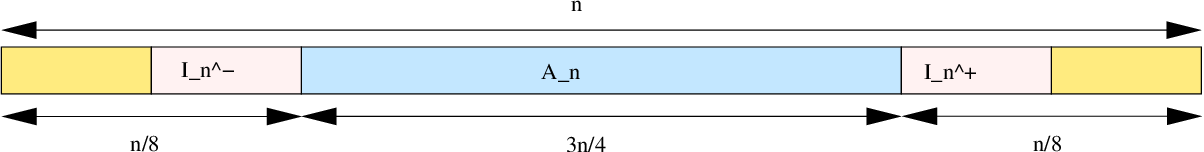}
\end{figure}

Given $\alpha_n^-\in I_n^-$ and $\alpha_n^+\in I_n^+$, define
cylinders
\begin{align*}
  B(\alpha_n^-)&=\left\{[x+iy]\in T_n\ \bigl|\bigr.\ y\leq\alpha_n^-\right\},\\
  B(\alpha_n^+)&=\left\{[x+iy]\in T_n\ \bigl|\bigr.\ \alpha_n^+\leq y\right\}.
\end{align*}

\begin{figure}[h]
  \centering
  \psfrag{B(a_n^-)}[][][.7]{$B(\alpha_n^-)$}
  \psfrag{B(a_n^+)}[][][.7]{$B(\alpha_n^+)$}
  \psfrag{a_n^-}[][][.7]{$\alpha_n^-$}
  \psfrag{a_n^+}[][][.7]{$\alpha_n^+$}
  \includegraphics[width=10cm]{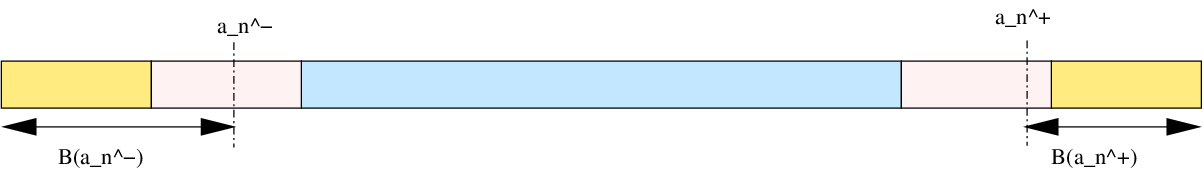}
\end{figure}

Their lengths are at least $n/16$. Define
$w_n:T_n\rightarrow\mathbb{R}$ by the following differential problem:
$$\begin{cases}
  \Delta w_n=0 &\mbox{ on } B(\alpha_n^-)\cup B(\alpha_n^+),\\
  w_n=0 &\mbox{ on } \partial T_n=\left\{[x+iy]\in\mathbb{C/Z}\ \Bigl|\Bigr.\ |y|=n/2\right\},\\
  w_n=u_n &\mbox{ on }T_n\setminus\left(B(\alpha_n^-)\cup B(\alpha_n^+)\right).
\end{cases}
$$
Since the continuous function $w_n$ is piecewise smooth and satisfies
$w_n=0$ on the boundary of $T_n$, it is compatible with the
identification of the boundary and induces a piecewise smooth function
on the torus.

Let $\delta_n=\int_{T_n}w_n\,dg_n$.
Since $\int_{A_n^\epsilon}w_n\,dg_n=0$, it follows from
concentration of the measures $dg_n$ on $A_n^\epsilon$ and from the maximum
principle that
\begin{gather*}
  |\delta_n|
  =\left|\int_{T_n\setminus
    A_n^\epsilon}w_n\,dg_n\right|
  \leq\max_{x\in T_n\setminus  A_n^\epsilon}|w_n(x)|\int_{T_n\setminus A_n^\epsilon}dg_n
  \leq \varepsilon.
\end{gather*}
This means that $w_n$ is almost admissible for the variational
characterization of $\lambda_1(g_n)$. 
Define $f_n:T_n\rightarrow\mathbb{R}$ by 
\begin{gather}\label{deftorusfn}
  f_n=w_n-\delta_n
\end{gather}
so that
$\int_{T_n}f_n\,dg_n=0$.
From \eqref{renormcentermass} and \eqref{tagada} it follows that
\begin{align}\label{bornedenominateur}
  \int_{T_n}f_n^2\,dg_n&\geq \int_{A_n^\epsilon}(u_n-\delta_n)^2\,dg_n\\
  &=\int_{A_n^\epsilon}u_n^2\,dg_n+\delta_n^2\int_{A_n^\epsilon}dg_n
  \geq \frac{\mathcal{A}(A_n^\epsilon,g_n)}{3}
  \geq \frac{1-\epsilon}{3}\nonumber
\end{align}
Using the variational
characterization~\eqref{varcharaclambda1} of $\lambda_1(g_n)$ and
conformal invariance of the Dirichlet energy leads to
\begin{align}\label{varlambdaineq}
  \lambda_1(g_n)&
  \leq \frac{3\int_{T_n}|\nabla f_n|^2\,dg_n}{1-\varepsilon}\\
  &=\frac{3}{1-\varepsilon}
  \left(\int_{T_n\setminus(B(\alpha^-_n)\cup B(\alpha^+_n))}|\nabla w_n|^2\,dg_n\right.\nonumber\\
    &\hspace{4.55cm}\left.+\int_{B(\alpha^-_n)\cup B(\alpha^+_n)}|\nabla w_n|^2\,dg_n\right)\nonumber\\
  &\leq
  \frac{3}{1-\varepsilon}
  \left(\int_{S^2}|\nabla\pi_i|^2\,dg_{S^2}
    +\int_{B(\alpha^-_n)\cup B(\alpha^+_n)}|\nabla w_n|^2\,dg_n\right)\nonumber\\
  &=
  \frac{8\pi}{1-\varepsilon}+\frac{3}{1-\varepsilon}
  \int_{B(\alpha^-_n)\cup B(\alpha^+_n)}|\nabla w_n|^2\,dg_n.\nonumber
\end{align}

\subsection{Energy estimate}
On a long flat cylinder like $B(\alpha_n^\pm)$, the Dirichlet energy
of a harmonic function is controled by the Dirichlet energy of its
restriction to the boundary circles.
Corollary~\ref{coro_nrjdirich} is applied to the function $f=w_n$
on the cylinders
$\Omega=B(\alpha^\pm_n)$, their lengths are at least $n/16$. For any
$x$, $u(x):=f(x,\alpha_n^\pm)\in[-1,1]$ since it is a coordinate
function on the sphere.

The next lemma shows that the numbers $\alpha_n^\pm\in I_n^\pm$ can be chosen to make the
Dirichlet energy of $u_n$ on the boundary of $B(\alpha^\pm_n)$ small.
\begin{lemma}\label{nrjboundary}
  There exist $\alpha_n^-\in I_n^-$ and $\alpha_n^+\in I_n^+$
  such that
  $$\int_{x=0}^{1}|\partial_xu_n(x+i\alpha_n^\pm)|^2dx\leq \frac{128\pi}{3n}.$$
\end{lemma}
\begin{proof}
  Let
  $E_n^\pm=\left\{[x+iy]\ \Bigl|\Bigr.\ y\in I_n^\pm\right\}.$
  Since the width of $I^\pm_n$ is at least~$n/16$, the mean-value
  theorem implies
  \begin{gather*}
    \iint_{E_n^\pm}|\nabla u_n|^2\ dxdy\geq\frac{n}{16}\min\left\{\int_{x=0}^{1}|\nabla u_n(x+i\alpha^\pm_n)|^2dx\
      \Bigl|\Bigr.\ \alpha_n^\pm\in I_n^\pm\right\}
  \end{gather*}
  In particular, there exists $\alpha^\pm_n\in I^\pm_n$ such that
  \begin{gather*}
    \int_{x=0}^{1}|\partial_xu_n(x+i\alpha_n^\pm)|^2dx\leq\frac{16}{n}\iint_{E_n^\pm}|\nabla u_n|^2\ dxdy.
  \end{gather*}
  By conformal invariance of the Dirichlet energy, this is bounded
  above by 
  \begin{gather*}
    \frac{16}{n}\int_{\tau_n\circ\phi(E^{\pm}_n)}|\nabla\pi_i|^2\,dg_{S^2}
    \leq\frac{16}{n}\int_{S^2}|\nabla \pi_i|^2\,dg_{S^2}=\frac{16}{n}.\frac{8\pi}{3}=\frac{128\pi}{3n}.
  \end{gather*}
\end{proof}
\begin{proof}[Proof of Theorem~\ref{monthmNad}]
  Let $f_n$ be the function given by~\eqref{deftorusfn}.
  Using the estimate on boundary derivative (Lemma~\ref{nrjboundary}),
  the Dirichlet energy estimate on cylinders (Lemma~\ref{nrjdirich})
  implies
  \begin{gather}\label{tubestimatefn}
    \lim_{n\rightarrow\infty}\int_{B(\alpha_n^\pm)}|\nabla f_n|^2\,dg_n=0.
  \end{gather}
  Using inequality~(\ref{varlambdaineq}), obtained by the variational
  characterization of $\lambda_1(g_n)$, it follows that
  $\limsup_{n}\lambda_1(g_n)\leq\frac{8\pi}{1-\epsilon}.$
  Since $\epsilon>0$ is arbitrary, it follows that
  $$\limsup_{n\rightarrow\infty}\lambda_1(g_n)\leq 8\pi.$$
  
  Finally, the lower bound  follows from the result of Friedlander and
  Nadirashvili \cite{nad:2} stated in the introduction: for any
  conformal class $C$ on a closed surface, $\nu(C)\geq 8\pi$.
\end{proof}

\section{Conformal degeneration on the Klein bottle}\label{sectionklein}
The goal of this section is to prove Theorem~\ref{kleindegeneration}.
Let $S^k(r)$ be the $k$-dimentional sphere of radius $r$ with its
standard metric $g_{S^k(r)}$ and
$\mathbb{R}P^2(r)$ be the associated projective plane with standard
metric $g_{\mathbb{R}P^2(r)}$.
Recall from the introduction that any Klein bottle is conformally
equivalent to a unique $K_b=\mathbb{C}/G_b$ where $G_b$ is the group of
transformations of $\mathbb{C}$ generated by $t_b(x+iy)=x+i(y+b)$ and
$\tau(x+iy)=x+\pi-iy$. The rectangle $[0,\pi]\times[-b/2,b/2]$ is a
fundamental domain for the action of $G_b$ on $\mathbb{C}$.
Reversing and identifying the opposite vertical sides this rectangle,
we obtain  a Möbius strip
$$M_b=\bigl([0,\pi]\times[-b/2,b/2]\bigr)/\tau.$$

\subsection{Transplantation to the sphere $S^4$ via projective space}
In this paragraph we exhibit a conformal embedding of the 
infinite Möbius strip $M_\infty$ in the sphere $S^4$. We start with a
lemma which will be used to embed Möbius strip conformally in
$\mathbb{R}P^2$.
\begin{lemma}
  The conformal application
  $\phi:\mathbb{C}\rightarrow
  S^2\subset\mathbb{R}^3$
  defined by
  $$\phi(x+iy)=\frac{1}{e^{2y}+1}(2e^y\cos(x), 2e^y\sin(x), e^{2y}-1)$$
  satisfies 
  $$\phi(x+2\pi+iy)=\phi(x+iy),$$
$$\phi(\tau(x+iy))=-\phi(x+iy).$$
  It induces a conformal equivalence
  $\phi:M_\infty\rightarrow\mathbb{R}P^2\setminus\left\{[0:0:1]\right\}$.
\end{lemma}

The Veronese map $v$ is a well known minimal isometric embedding of
$\mathbb{R}P^2(\sqrt{3})$ in the sphere $S^4$ by its first
eigenfunctions. This means that the components of $v$ are eigenfunctions for
$\lambda_1(g_{\mathbb{R}P^2(\sqrt{3})})=2.$ For details, see
\cite{berger:1} and \cite{lawson:1}.

It follows that the composition $v\circ\phi$ is a conformal
embedding of the Möbius strip $M_\infty$ in $S^4$.

\subsection{Concentration on Möbius strips}
Without loss of generality, consider a sequence $b_n=n$.
\begin{lemma}\label{mobconcentration}
  If there exists $K>0$ such that
  $$\liminf_{n\rightarrow\infty}\lambda_1(g_n)\geq K,$$
  then for any $\epsilon>0$, there exists $N\in\mathbb{N}$ such that for
  any $n\geq N$, 
  $$\max\bigl\{\int_{M_{3n/4}}dg_n,\int_{K_n\setminus K_{3n/4}}dg_n\bigr\}\geq 1-\epsilon.$$
\end{lemma}
\begin{proof}
  Since the function
  $$u_n([x+iy])=\cos(\frac{2\pi y}{n})$$
  used in the proof of Lemma~\ref{lemmeConcentration} is even, it induces
  a first eigenfunction on the flat Klein bottle $K_n$.
  The cylinders constructed in this proof are compatible with the
  identification on the Möbius strip, also because the function $u_n$
  is even.
\end{proof}

\noindent
\textbf{Notation}: Without loss of generality, we suppose the maximum
is reached by $A_n^\epsilon:=M_{3n/4}.$ See
Lemma~\ref{lemmeConcentration} for details.

\subsection{Renormalization of the centers of mass}

It follows from Corollary~\ref{coroAlaHerschplongement} that there exists
conformal transformations $\tau_n$ of $S^4$ such that
$$\int_{A_n^\epsilon}\pi\circ\tau_n\circ v\circ\phi\,dg_n=0$$
where $\pi:S^4\hookrightarrow\mathbb{R}^5$ is the standard
embedding. 
For $1\leq i\leq 5$, let
\begin{gather}\label{defkleinuni}
  u^i_n=\pi_i\circ\tau_n\circ v\circ\phi.
\end{gather}

\subsection{Test functions}
The numbers
$$\alpha_n\in I_n:=[\frac{6n}{16},\frac{7n}{16}],$$
will be chosen later (see Lemma~\ref{nrjboundary}).
For each $1\leq i\leq 5$, define $w^i_n:M_n\rightarrow\mathbb{R}$
by the following differential problem:
$$\begin{cases}
  \Delta w^i_n=0 &\mbox{ on } M_{n}\setminus M_{\alpha_n},\\
  w^i_n=0 &\mbox{ on } \partial M_{n}=\left\{[(x,y)]\in M_{n}\ \bigl|\bigr.\ |y|=n/2\right\},\\
  w^i_n= u^i_n&\mbox{ on }M_{\alpha_n}
\end{cases}
$$
where $u^i_n$ is defined by \eqref{defkleinuni}.
Since the continuous function $w^i_n$ is piecewise smooth and satisfies
$w^i_n=0$ on the boundary of $M_n$, it is compatible with the
identification of the boundary and induces a piecewise smooth function
on the Klein bottle $K_n$.

Since
$\int_{A_n^\epsilon}w^i_n\,dg_n=0,$
it follows from
concentration of the measures $dg_n$ on $A_n^\epsilon$ and from the maximum
principle that $\delta^i_n:=\int_{M_{n}}w^i_n\,dg_n$
satisfies $|\delta^i_n|\leq\epsilon$.
This means that $w^i_n$ is almost admissible for the variational
characterization of $\lambda_1(g_n)$. Thus, it is natural
to define $f^i_n:M_{n}\rightarrow\mathbb{R}$ by
\begin{gather}\label{defkleinfni}
  f^i_n=w^i_n-\delta^i_n
\end{gather}
so that for each $i$
$$\int_{M_{n}}f^i_n\,dg_n=0$$
and similarly to inequality~(\ref{bornedenominateur})
\begin{gather}\label{inegalitesommefn}
  \sum_{i=1}^5\int_{M_n}(f^i_n)^2\,dg_n\geq\int_{A_n^\epsilon}dg_n\geq 1-\epsilon.
\end{gather}

It follows from the variational characterization of $\lambda_1(g_n)$
and from inequality~(\ref{inegalitesommefn}) that
\begin{align}\label{varcarklein}
  \lambda_1(g_n)\left(1-\epsilon\right)
  &\leq
  \lambda_1(g_n)\left(\sum_{i=1}^5\int_{M_n}(f^i_n)^2\,dg_n\right)
  \leq\sum_{i=1}^5\int_{M_n}|\nabla f^i_n|^2\,dg_n\nonumber\\
  &=\sum_{i=1}^5\left(
  \int_{M_{\alpha_n}}|\nabla w^i_n|^2\,dg_n+
  \int_{M_n\setminus M_{\alpha_n}}|\nabla w^i_n|^2\,dg_n\right).
\end{align}

\subsection{Energy estimate}
\begin{itemize}
\item[]
\end{itemize}
\noindent
{\em First step: Bounding
$\sum_{i=1}^5\int_{M_{\alpha_n}}|\nabla w^i_n|^2\,dg_n.$}\\

Recall that the Veronese embedding
$v:\mathbb{R}P^2(\sqrt{3})\rightarrow S^4$ is isometric and minimal.
On $\Sigma_n:=\tau_n\circ v(\mathbb{R}P^2(\sqrt{3}))$ we consider the
metric induced by $g_{S^4}$.
Proposition~1 of~\cite{li:1} says that if a compact surface is
minimally immersed in a sphere, then its area cannot be increased by
conformal transformations of the sphere. In our particular case this
leads to the following proposition.
\begin{proposition}\label{LiYauminsurface}
  For each $n\in\mathbb{N}$,
  \begin{gather*}
    \int_{\Sigma_n}dg_{S^4}\leq \mbox{Area of }\mathbb{R}P^2(\sqrt{3})=6\pi.
  \end{gather*}
\end{proposition}
It follows from conformal invariance of the Dirichlet energy that
\begin{align*}
  \sum_{i=1}^5
  \int_{M_{\alpha_n}}|\nabla w^i_n|^2\,dg_n
  &\leq \sum_{i=1}^5
  \int_{\Sigma_n}|\nabla \pi_i|^2\,dg_{S^4}
\end{align*}
It is proved on page 146 of~\cite{schoen:1} that
$\Sigma_n$ being isometrically immersed in $S^4$
implies the point-wise identity $\sum_{i=1}^5|\nabla \pi_i|^2=2$.
Whence
\begin{gather}\label{kleinEnergieEst1}
  \sum_{i=1}^5\int_{M_{\alpha_n}}|\nabla w^i_n|^2\,dg_n\leq2\int_{\Sigma_n}dg_{S^4}\leq 12\pi.
\end{gather}

\noindent
{\em Second step: Bounding
$\sum_{i=1}^5\int_{M_n\setminus M_{\alpha_n}}|\nabla w^i_n|^2\,dg_n.$}\\

The function $w_n$ is harmonic on the set $M_{n}\setminus
M_{\alpha_n}$. This is a cylinder of length
$L_n:=n-\alpha_n\geq 9/16n$ and of width $2\pi$.
The next lemma shows that the Dirichlet energy of $w_n$ on the
boundary of these cylinders can be controlled by appropriate choice of
$\alpha_n$.
\begin{lemma}\label{nrjboundarymob}
  The number $\alpha_n\in I_n=[6n/16,7n/16]$
  can be chosen such that
  $$\sum_{i=1}^5\int_{x=0}^{\pi}|\partial_xw^i_n(x\pm i\alpha_n)|^2dx\leq 192\pi/n.$$
\end{lemma}
\begin{proof}
  We argue as in  Lemma~\ref{nrjboundary}.
  Let $E_n^\pm=\{[x+iy]\ \left|\right.\ y\in I_n\}.$
  It follows from the mean-value theorem that
  \begin{align*}
    \frac{n}{16}\min_{\alpha_n^\pm\in I_n^\pm}
      \sum_{i=1}^5\int_{x=0}^{2\pi}|\nabla w^i_n(x,\alpha^\pm_n)|^2dx
    &\leq\sum_{i=1}^5\iint_{E_n^\pm}|\nabla w^i_n(x+iy)|^2\ dxdy\\
    &\leq 12\pi.
  \end{align*}
\end{proof}
Since $M_n\setminus M_{\alpha_n}$ has length $L_n$ at least $n/16$ and is of
width $2\pi$, Lemma~\ref{nrjdirich} implies
\begin{align}\label{kleinEnergieEst2}
  \sum_{i=1}^5\int_{M_n\setminus M_{\alpha_n}}|\nabla w^i_n|^2\,dg_n&\\
  &\hspace{-3cm}\leq\sum_{i=1}^5\Bigl(\frac{2\pi}{L_n}+\coth(L_n)\int_{x=0}^{2\pi}|\partial_xw^i_n(x\pm i\alpha_n)|^2\,dx\Bigr)\nonumber\\ 
  &\leq\frac{\pi}{n}\left(160+5\times 192\coth(\frac{n}{16})\right).\nonumber
\end{align}

\begin{proof}[Proof of Theorem~\ref{kleindegeneration}]
  Substitution of inequality~\eqref{kleinEnergieEst1} and
  inequality~\eqref{kleinEnergieEst2} in the variational
  characterization~\eqref{varcarklein} leads to
  $$\limsup_{n\rightarrow\infty}
  \lambda_1(g_n)(1-\epsilon)\leq 12\pi+\limsup_{n\rightarrow\infty}\frac{\pi}{n}\left(160+960\coth(n/16)\right)=12\pi.$$
  Since $\epsilon>0$ is arbitrary, this completes the proof of
  Theorem~\ref{kleindegeneration}.
\end{proof}
\section{Concentration to points}\label{sectionPunkt}
The main goal of this section is to prove Theorem~\ref{punktthm}.
We start by proving that concentration to a point has
no influence on spectrum.
\begin{proof}[Proof of Proposition~\ref{no}]
  There exists a sequence of diffeomorphisms $\phi_n$ such that
  $\lim_{n\rightarrow\infty}\phi_n(x)=p\ \ dg_0$-almost everywhere. Indeed, 
  let $f:M\rightarrow\mathbb{R}$ be any Morse function having $p$ as
  its unique local minimum and consider
  $\phi\mathbb{R}\times M\rightarrow M$ be its negative gradient
  flow with respect to $g_0$. Since the stable  manifolds of any
  critical point other than $p$ are of codimension strictly 
  greater than 1, they are $dg_0$-negligible, hence the sequence
  $\phi_n$ satisfies the required property and
  $g_n=\phi_n^*g_0$ concentrates to $p$.
\end{proof}

We know proceed with the proof of Theorem~\ref{punktthm}.
\subsection{Construction of a neighborhood system}
Let $(g_n)\subset\mathcal{R}(\Sigma)$ be a
sequence of metrics concentrating to $p\in\Sigma$.
Let $\mathbb{D}$ be the unit open disk in $\mathbb{C}$.
Let $\eta:D\rightarrow\mathbb{D}$ a conformal chart
arround $p$ such that $\eta(p)=0$. Observe that since the metrics
$g_n$ are all in the same conformal class, the same chart $\eta$ will
be conformal for each of them.
\begin{lemma}
  There exists a conformal equivalence
  $$\psi:D\setminus\{p\}\rightarrow(0,\infty)\times S^1$$
  and a family $\mathcal{U}_n\subset D$ of neighborhood of $p$ such that
  $$\lim_{n\rightarrow\infty}\int_{\mathcal{U}_n}dg_n=1$$
  and
  $$\psi(D\setminus\mathcal{U}_n)=(0,L_n)\times S^1$$
  with $L_n\rightarrow\infty$ and $\psi(\partial D)=\{0\}\times S^1.$
\end{lemma}
\begin{proof}
For $0\leq\epsilon\leq1$, let  $\mathcal{U}(\epsilon)=\eta^{-1}B(0,\epsilon)$.
Define $\epsilon_n$ by
$$\int_{\mathcal{U}(\epsilon_n)}dg_n=1-\epsilon_n.$$
It follows from concentration that
$\lim_{n\rightarrow\infty}\epsilon_n=0$
so that
$$\lim_{n\rightarrow\infty}\int_{\mathcal{U}(\epsilon_n)}dg_n=1$$
Define the conformal equivalence
$\alpha:\mathbb{R}\times S^1\rightarrow\mathbb{C}^*$ by
$$\alpha(x,y)=e^{-x}\bigl(\cos(y),\sin(y)\bigr).$$
The composition $\psi=\alpha^{-1}\circ\eta$ 
has the required property since
$$\psi(D\setminus\mathcal{U}(\epsilon_n))=\bigl(0,-\ln(\epsilon_n)\bigr)\times S^1.$$  
\end{proof}

\subsection{Renormalization of the centers of mass}
Recall that $\sigma$, as defined in \eqref{stereoparam}, is the
stereographic parameterization of the sphere $S^2$ by its equatorial
plane. Let $H\subset S^2$ be the southern hemisphere of $S^2$.
The map
$$\phi=\sigma\circ\eta:D\rightarrow H$$
is a conformal equivalence such that $\phi(p)$ is the south pole.
It follows from Corollary~\ref{coroAlaHerschplongement} that
there exist conformal transformations $\tau_n$ of the sphere such that
$$\int_{\mathcal{U}_n}\pi\circ\tau_n\circ\phi\,dg_n=0$$
where $\pi:S^2\hookrightarrow\mathbb{R}^3$ is the standard embedding.
For each $n\in\mathbb{N}$,
since $\pi_1^2+\pi_2^2+\pi_3^2=1$ on $S^2$, there exists an indice $i=i(n)\in\{1,2,3\}$
such that the function
$u_n=\pi_i\circ\tau_n\circ\phi$ satisfies
$$\int_{\mathcal{U}_n}u_n^2\,dg_n\geq\frac{1}{3}\int_{\mathcal{U}_n}dg_n.$$

\subsection{Test functions}
Consider $\alpha_n\in[L_n/2,L_n]$ to be chosen later.
Define $w_n:\Sigma\rightarrow\mathbb{R}$ as the unique solution of
$$\begin{cases}
  w_n=0&\mbox{ on }\Sigma\setminus D\\
  w_n=\pi_i\circ\tau_n\circ\phi&\mbox{ on } 
  \mathcal{U}_n\cup\psi^{-1}\bigl((\alpha_n,L_n)\times S^1\bigr)\\
  \Delta w_n=0&\mbox{ on } \psi^{-1}\bigl((0,\alpha_n)\times S^1\bigr).
\end{cases}$$
By the maximum principle, $\delta_n:=\int_\Sigma w_n\,dg_n\leq
\int_{\Sigma\setminus\mathcal{U}_n}dg_n$.
Define $f_n=w_n-\delta_n$ so that $\int_{\Sigma}f_n=0$ and  $f_n$ is admissible
for the variational characterisation of $\lambda_1(g_n)$:
\begin{align*}
  \lambda_1(g_n)&\leq\frac{\int_{\Sigma}|\nabla
    f_n|^2\,dg_n}{\int_{\Sigma}f_n^2\,dg_n}
  \leq\frac{\int_{D}|\nabla w_n|^2\,dg_{S^2}}
  {\int_{\Sigma}(w_n-\delta_n)^2\,dg_n}.
\end{align*}
The denominator satisfies
\begin{align*}
  \int_{\Sigma}(w_n-\delta_n)^2\,dg_n&\geq\int_{\mathcal{U}_n}(w_n^2-2\delta_nw_n+\delta_n^2)\,dg_n\\
  &\geq\int_{\mathcal{U}_n}w_n^2\,dg_n\\
  &\geq\frac{1}{3}\int_{\mathcal{U}_n}dg_n.
\end{align*}
Whence,
\begin{align*}
  \frac{\lambda_1(g_n)}{3}\int_{\mathcal{U}_n}dg_n&\leq\int_{D}|\nabla w_n|^2\,dg_n\\
  &\leq\int_{\mathcal{U}_n\cup\psi_n^{-1}(0,\alpha_n)\times S^1}|\nabla w_n|^2\,dg_n+
  \int_{\psi_n^{-1}(\alpha_n,L_n)\times S^1}|\nabla w_n|^2\,dg_n\\
  &\leq\frac{8\pi}{3}+
  \int_{\psi_n^{-1}(\alpha_n,L_n)}|\nabla w_n|^2\,dg_n.
\end{align*}

\begin{proof}[Proof of Theorem~\ref{punktthm}]
  The set $\psi^{-1}\bigl((0,\alpha_n)\times S^1\bigr)$ where $w_n$ is
  harmonic is conformally equivalent to a cylinder of length
  $\alpha_n\geq\frac{L_n}{2}$ which becomes infinite as $n$ goes to
  infinity. The proof is completed by choosing appropriate $\alpha_n$
  as it was done in Lemma~\ref{nrjboundary} and
  Lemma~\ref{nrjboundarymob} and then applying Lemma~\ref{nrjdirich}
  to bound the Dirichlet energy.
\end{proof}

\section{Acknowledgment}
This work is a part of my PhD thesis at the Université de Montréal
under the supervision of Marlène Frigon and Iosif Polterovich. I want
to thank them for their support and encouragement. I would like to
thank N. Nadirashvili for outlining the proof of 
Theorem~\ref{monthmNad} and~\ref{kleindegeneration} and for
interesting discussions. I would also like to thank A. El Soufi,
P. Guan, S. Ilias, M. Levitin, and M. Sodin for interesting and useful
conversations. Finally, I want to thank my friends Baptiste Chantraine, Nicolas
Beauchemin, Rémi Leclercq, and Jean-François Renaud for useful
comments on an earlier version of this paper.

\end{document}